\documentclass[12pt]{article}

\title{On Gorenstein surfaces dominated by ${\bf P}^2$}
\author{R.~V.~Gurjar,~ C.~R.~Pradeep,~ D.~-Q.~Zhang}
\date{}

\newcommand{\ncom}{\newcommand}
\ncom{\se}{\setcounter{equation}{0}}

\ncom{\nind}{\noindent}
\ncom{\comx}{{\bf C}}
\ncom{\be}{\begin{equation}}
\ncom{\ee}{\end{equation}}
\ncom{\beqn}{\begin{eqnarray*}}
\ncom{\eeqn}{\end{eqnarray*}}
\ncom{\beq}{\begin{eqnarray}}
\ncom{\eeq}{\end{eqnarray}}
\ncom{\nno}{\nonumber}
\ncom{\hs}{\mbox{\hspace{.25cm}}}
\ncom{\rar}{\rightarrow}
\ncom{\Rar}{\Rightarrow}
\ncom{\noin}{\noindent}
\ncom{\bc}{\begin{center}}
\ncom{\ec}{\end{center}}
\ncom{\zee}{{\bf Z}}
\ncom{\aff}{\mathbb A}
\ncom{\proj}{{\bf P}}
\ncom{\sz}{\scriptsize}
\ncom{\f}{\frac}
\ncom{\la}{\lambda}
\ncom{\ka}{\overline{\kappa}}
\ncom{\Ga}{\Gamma}
\ncom{\ga}{\gamma}
\ncom{\si}{\sigma}
\ncom{\Si}{\Sigma}
\ncom{\bib}{\bibitem}
\ncom{\sst}{\subset}
\ncom{\sms}{\setminus}
\ncom{\est}{\emptyset}
\ncom{\pf}{\noin Proof. }
\ncom{\bighs}{\hspace{.5 cm}}
\ncom{\ulin}{\underline}
\ncom{\olin}{\overline}
\ncom{\dm}{\displaystyle}
\ncom{\cstar}{{\bf C}^*}
\ncom{\vtld}{\tilde{V}}
\ncom{\ztld}{\tilde{Z}}
\ncom{\idealm}{\mathfrak m}
\ncom{\vnot}{V^\circ}
\ncom{\tnot}{T^\circ}
\ncom{\ctld}{\tilde{C}}
\ncom{\De}{\Delta}
\ncom{\om}{\omega}

\newtheorem{theorem}{Theorem}
\newcommand{\bt}{\begin{theorem}}
\newcommand{\et}{\end{theorem}}

\begin{document}
\maketitle

\begin{small}

\parskip = 10pt

\nind
{\bf Abstract.}  In this paper we prove that 
a normal Gorenstein surface dominated by $\proj^2$ is 
isomorphic to a quotient $\proj^2/G$, 
where $G$ is a finite group of automorphisms of $\proj^2$ (except possibly 
for one surface $V_8'$). We can completely classify all such quotients. 
Some natural conjectures when the surface is not Gorenstein are also stated.\\

\nind
{\it Mathematics Subject Classification (2000):} 14J25, 14J26, 14J45.\\

\nind
{\it Key words:} log del Pezzo surface, Gorenstein surface.

\section{Introduction}
Let $V$ be a normal projective surface defined over $\comx$. $V$ is said
to be a {\em log del Pezzo} surface if $V$ has at worst quotient
singularities and the anti-canonical divisor $-K_V$ is ample. The {\em
rank} of $V$ is the Picard number 
$\rho(V) = dim_{{\bf Q}} Pic(V) \otimes {\bf Q}$. It
is easy to see that any quotient of $\proj^2$ by a finite group of
automorphisms is a log del Pezzo surface of rank one. Miyanishi and Zhang
have raised the question of giving a criterion for a projective normal
surface to be isomorphic to $\proj^2/G$, where $G$ is a finite group of
automorphisms of $\proj^2.$ In [9] certain rank $1$ log del
Pezzo surfaces are shown to be quotients of $\proj^2$ modulo a finite
group. \\
Our main theorem is the following result.

\noin {\bf Theorem 1.} {\em Let $V$ be a Gorenstein normal surface and let
$f: \proj^2 \rar V$ be a non-constant morphism. Then we have the
following assertions: 

\noin (1) If $\pi_1(V \sms Sing~V)$ is trivial then $V$ is isomorphic
to one of the following surfaces. \\ 
(i) The projective plane $\proj^2$, \\ 
(ii) The quadric cone $Q := \{X^2 + Y^2 + Z^2 = 0 \}$ in
$\proj^3$, \\
(iii) A surface of singularity type $A_1+A_2$, or\\
(iv) The surface $V_8'$ which has a unique singular point, which is 
analytically the $E_8$ singular point. (cf. Theorem 2 below).\\
(2) If $\pi_1(V \sms Sing~V)$ is non-trivial then $f$ factors as
$\proj^2 \rar W \rar V$, where $W$ is one of the surfaces (i) or
(ii) in (1) above and $W \rar V$ is \'{e}tale over $V \sms Sing~V$.\\
(3) $V$ is always isomorphic to a quotient $\proj^2/G$ for a
finite group of automorphisms of $\proj^2$, except for the surface 
$V_8'$ in the case (iv) above. The surface $V_8'$ is not isomorphic to 
any quotient of $\proj^2$ modulo a finite group of automorphisms.}\\

\noin {\bf Remarks.} (1) We can give a very precise description of any $V$
in part (3) above, particularly its singularity type and 
$\pi_1(V \sms Sing~V)$ and the corresponding surface $W$ as in part
(2) above. (See, \S  5.)\\
(2) The surfaces $V_8,~V_8'$ are ``twin'' surfaces. Theorem 1 says that 
there is no non-constant morphism $\proj^2\rar V_8.$ It is most probable 
that there is no such map $\proj^2\rar V_8'$ but we have been unable to 
prove this. This is the exception in p
art (iv) above.\\

We will also prove the following result which will be used in the proof of
Theorem 1. This result is stated in [9] and a sketch of proof is given
there. In view of the importance of this result for our proof of Theorem 1
we will give a complete proof. We would like to point out that the
uniqueness assertion made in Lemma 7 of [9] is not quite correct. For
$K_V^2=1$ we have found two non-isomorphic surfaces.

\noin {\bf Theorem 2.} {\em Let $V$ be a Gorenstein log del Pezzo
surface of rank 1 such that $\pi_1(V \sms Sing~V) = (1)$. Let
$d:=(K_V)^2$. Then we have the following assertions:

(i) $1 \leq d \leq 9$ and $d = 9$ implies $V \cong
\proj^2$. If $d = 1$ then $V$ is one of the following two
surfaces in the weighted projective space $\proj(1,1,2,3)$.
$$ V_8:  \{W^2 + Z^3 + X^5Y = 0 \}.$$
$$ V_8': \{W^2 + Z^3 + X^5Y + X^4Z = 0 \}.$$
Both these have a unique singularity of type $E_8$. The surface $V_8$ 
contains a rational curve $C$ with only one ordinary cusp (and
otherwise smooth) such that 
$C \sim -K_V$. The surface $V_8'$ contains a rational curve $C$
with only one ordinary double point (and otherwise smooth) such that 
$C\sim -K_V.$ In both these cases $C$ does not pass through the
singular point of $V.$

(ii) If $d > 1$, then $V$ contains a cuspidal rational curve $C$
as in (i) above such that $C \sim -K_V$  and $C$ does not pass through the
singular
point of $V$. The surface $V$ is uniquely determined by the integer $d$. 

(iii) If $d\neq 8$ then $V$ contains an irreducible curve $\De$ such that 
$\De$ generates the Weil divisor class group of $V$ and $-K_V \sim d \De$.
If $d=8$ then $V$ contains an irreducible curve $\De$ which generates the 
Weil divisor class group of $V$ and $-K_V\sim 4\De$.

(iv) In case $1<d\leq 5$ or $d=1$ and $V=V_8$ the affine surface 
$V\sms C$ is isomorphic to $\comx^2/G,$ where $G$ is a 
finite group of automorphisms of
$\comx^2$ isomorphic to the fundamental group at infinity of
$V \sms C$.}\\

We now mention several results proved by other mathematicians which are
closely related to our Theorems 1 and 2.\\
(1) Demazure has proved important results about Gorenstein log del
Pezzo surfaces in [3], particularly about the linear systems $|-nK_V|.$\\
(2) In [4] general results about embeddings of Gorenstein log del
Pezzo surfaces
are proved.\\
(3) In [9] a classification of rank 1 Gorenstein log del Pezzo
surfaces is given.\\
(4) In [7] R. Lazarsfeld has proved that any smooth variety which is
dominated by $\proj^n$ is isomorphic to $\proj^n.$\\
(5) In [5] and [8] it is proved that if there is a proper map
$f:\comx^2\rar V$ onto a normal algebraic surface $V$ then $V$ is
isomorphic to $\comx^2/G$ for a finite group of automorphisms of
$\comx^2.$\\ 
(6) In [6], Mohan Kumar has shown that 
if a normal rational surface $S$ has a singularity of type $E_8$ 
then the local ring of $S$ at this singular point is 
isomorphic to one of the (non-isomorphic) local rings,
viz. $\comx~[X,Z,W]_{(X,Z,W)}/(W^2+Z^3+X^5)$ or, 
$\comx~[X,Z,W]_{(X,Z,W)}/(W^2+Z^3+X^4Z+X^5).$ 
This result will be used in {\S 3}.\\
(7) In [1] it is proved that a surface of the form 
$\proj^2/G$ cannot have a unique singular point of 
the type $E_6,E_7$ or $E_8.$ This is a special case of our Theorem 1.\\

Our proof of Theorem 1 is almost self-contained. We do not use the
classification mentioned in (3) above. All we need is some general results
about embeddings given by $|-K_V|$ and $|-2K_V|$ which are proved in
[4]. Nevertheless, the paper [9] 
has been important for us while 
thinking of the proofs in this paper.\\
Recently, the first named author has proved the following general result
using the results of this paper in an important way.\\
{\it Let $\pi:{\proj}^2\rar{\proj}^2$ be a non-constant morphism. Let $C$
be
an irreducible curve of degree $>1$ in ${\proj}^2$ which is ramified for
$\pi$. Then the greatest common divisor of the ramification indices of the
irreducible curves lying over $C$ is $1$. In particular, $\pi^{-1}(C)$
cannot be irreducible.}\\

\section{Preliminaries}

All the algebraic varieties we consider are defined 
over the field $\comx$ of complex numbers.

A smooth complete rational curve $C$ on a smooth algebraic surface $S$ is
called a $(-n)$-{\it curve} if $C^2=-n.$ 
 
Let $Z$ be an irreducible normal variety such that $\pi_1(Z\sms sing~Z)$
is finite. Let $W'$ be the universal covering space of $Z\sms Sing~Z.$
Then $W'$ is also a variety. The normalization, $W$, of $Z$ in the
function field of $W'$ is called the {\it quasi-universal cover} of $Z.$ 
There is a proper morphism with finite fibers $W\rar Z$ 
which is unramified over $Z\sms Sing~Z.$

\nind
For any normal variety $Z$ we will denote the Zariski-open 
subset $Z\sms Sing~Z$ by $Z^0.$ 

\nind 
Let $Z$ be an irreducible normal variety. 
An algebraic action of $\cstar$ on $Z$ is said to be a {\it good} 
$\cstar$-action if $Z$ contains a point $z$ 
which is in the closure of every orbit. 
In this case we also say that $Z$ is a quasi-homogeneous variety.

In this section we collect a few results which we will use (at least
implicitly) to prove Theorems 1
and 2. The following is from lemma 6 of [9]. 
Recall that $d:=K_V^2.$ The integer $d$ is called the {\it degree} of $V$.

\noin {\bf Lemma 1.} (Reproved in Theorem 2.) 
{\em Let $V$ be a Gorenstein log del Pezzo surface
of rank one and let $\vnot$ denote the smooth locus of $V$. If
$\pi_1(\vnot)$ is trivial then $V$ has one of the following
combinations of singularities: $A_1, \hs A_1+A_2, \hs A_4, \hs D_5, 
\hs E_6, \hs E_7, \hs E_8$. The values of $K_V^2$ in these cases are
$8,~6,~5,~4,~3,~2,~1$ respectively. (Note that the case $K_V^2=7$ does
not occur.)}

Here $A_1+A_2$ means that there are two singular points on $V$, one of
type $A_1$ and the other of type $A_2$. We will call these as the
Dynkin types of $V$. Sometimes we also say that $V$ is of type $A_1,A_2.$

The following result plays an important role in our proof (see
Corollary 4.5 of [4]): 

\noin {\bf Lemma 2.} {\em
Let $V$ be a Gorenstein log del Pezzo surface of degree $d$.\\
If $d\geq 3$ then the linear system $|-K_V|$ is very ample and gives a
projectively normal embeding of $V$ in $\proj^d.$\\ 
If $d=2$ then $|-2K_V|$ is very ample and gives a projectively normal
embedding of $V$ in $\proj^6.$}

We will need the following result in our proof (see Lemma 6.1 of
[11]).

\noin {\bf Lemma 3.} {\em
Let $G$ be an algebraic group which acts algebraically on a normal variety
$Y.$ Suppose $f:Z\rightarrow Y$ is a finite morphism with $Z$ a normal
variety. If $Z$ contains a non-empty $G$-equivariant Zariski-open subset
$Z_0$ such that $f$ restricted to $Z_0$ is $G$-equivariant then
there is a unique action of $G$ on $Z$ such that $f$ is a $G$-morphism.}

\section{Proof of Theorem 2} 

In this section we will prove Theorem 2 (cf. Introduction). 
This and some of the arguments in the proof of Theorem 2 
will be used in proving Theorem 1.\\   
So let $V$ be a Gorenstein log del Pezzo surface of rank 1 such that
$\pi_1(V\sms Sing~V)$ is trivial. 

Let $g:\vtld\rar V$ be a minimal resolution of singularities. Then
$K_{\vtld}^2=K_V^2>0.$ By assumption, $|nK_{\vtld}|=\phi$ for $n>0.$
Since $\pi_1(V\sms Sing~V)$ is trivial the surface $\vtld$ is
simply-connected. Now it follows by Noether's theorem that $\vtld$
is rational. Hence $K_{\vtld}^2+~b_2(\vtld)=~10.$ This implies that $1\leq
K_{\vtld}^2\leq 9.$ If $b_2(\vtld)=1$ then $V$ is smooth and 
hence isomorphic to $\proj^2$. From now 
onwards, we will assume that $1\leq K_{\vtld}^2\leq 8.$\\
We will assume in what follows that $V$ is not smooth.

\nind
Consider first the case $d=8.$\\
Then $\vtld$ contains a $(-2)$-curve $C.$ This implies that $\vtld$ is the
Hirzebruch surface $\Sigma_2$ and $V$ is obtained by contracting the
$(-2)$-curve to an $A_1$-singularity. In this case $V$ is isomorphic to
the quadric cone $Q:=\{X^2+Y^2+Z^2=0\}$ in $\proj^3$.  
Let $x,y,z$ denote suitable
homogeneous coordinates on $\proj^2.$ Then the group $G:={\bf Z}/(2)$ acts
on $\proj^2$ by sending $[x,y,z]\rar [-x,-y,z].$ The line at infinity
$\{z=0\}$ is pointwise fixed and we see easily that $\proj^2/(G)$ is
isomorphic to $Q.$ Then Theorem 2 (ii), (iii) are clear.\\

\nind
{\bf Claim.} The case $d=7$ cannot occur.\\
To see this, assume that $d=7.$ Then $\vtld$ is obtained from a Hirzebruch
surface $\Sigma_n$ by one blowing-up with $E$ the exceptional curve. 
Let $S$ be the unique curve with
self-intersection $-n\leq 0$ on $\Sigma_n$ and $L$ be a fiber of the
$\proj^1$-fibration on $\Sigma_n.$ We can assume that the blown-up point
$p\in L.$ We have the formula $K_{\Sigma_n}\sim -2S-(2+n)L.$ If $p\not\in
S$ then $K_{\vtld}\sim -2S-(2+n)L'-(1+n)E,$ where $L'$ is
the proper transform of $L$. Clearly, $\vtld$
contains a $(-2)$-curve different from $S$, say $C.$ 
Since $K\cdot C=0,$ the curve $C$
is disjoint from $S,L,E.$ This is impossible.\\
If $p\in S$ then $K\sim -2S-(2+n)L'-(2+n)E.$ Again there is a $(-2)$-curve
$C$ different from $S.$ We get a contradiction as above.\\

In view of these observations, for the rest of the section we assume that
$1\leq d\leq 6.$\\ 

First we will give a construction of such surfaces with $1\leq d\leq 6$ and
later on prove that these are all the surfaces we are looking for.\\
 
We will construct examples of rank 1, Gorenstein log del Pezzo surfaces 
$V_i$, $3\le i \le 8$, and a rank 1, Gorenstein log del Pezzo surface 
$V_8'$ such that $K_{V_i}^2=9-i$, $K_{V_8'}^2=1$ ,
and of singularity types $A_1 + A_2, A_4, D_5, E_6, E_7, E_8, E_8$
respectively. Moreover, $V_i, V_8'$ are compactification 
of ${\bf C}^2$ with an irreducible boundary curve. 
In particular, these surfaces have simply-connected
smooth parts.

Let $X$ (resp. $X'$) be a relatively minimal rational elliptic surface
with a unique section $E$ ( see Claim 2 in Lemma 4) and
singular fibres of types $II^*, II$ (resp. $II^*, I_1, I_1$).
Such $X$ (resp. $X'$) is unique modulo fibration-preserving
isomorphisms. The construction and uniqueness is shown by
letting $X \rightarrow \proj^2$ (resp. $X' \rightarrow \proj^2$)
be the composition of blow-downs of the section $E$
and all components in the type $II^*$ fibre except for
a multiplicity $3$ component $C_3'.$ We get a pencil in $\proj^2$
generated by a cuspidal cubic and three times the tangent line at an
inflexion point (resp. a pencil generated by a nodal cubic and three
times the tangent line at an inflexion point). The pair of a cuspidal 
(resp. nodal) cubic
curve and the tangent line at an inflexion point is unique upto
projective transformations.\\
Write the type $II^*$ fibre as 
$\sum_{i=1}^6 iC_i + 4C_4' + 2C_2' + 3C_3'$ so that
$\sum_{i=1}^6 C_i + C_4' + C_2'$ is an ordered linear chain.

\par
Let $X \rightarrow V_8$ (resp. $X' \rightarrow V_8'$) be the contraction
of $E$ and the type $E_8$ divisor in the type $II^*$ fibre
to a smooth point and a singularity of type $E_8$.
Then $V_8$ and $V_8'$ are two non-isomorphic rank 1,
Gorenstein log del Pezzo surfaces of singularity type $E_8$.
Let $\De$ be the image on $V_8$ or $V_8'$ of $C_1$.

\par
For $3 \le i \le 8$ (resp. $i = 1$), we let $X \rightarrow V_i$ be the
contraction of $E$ and all components in the type $II^*$
fibre, except $C_{9-i}$ (resp. $C_4'$); denote by 
$\De$ the image on $V_i$ of $C_{9-i}$ (resp. $C_4$).

\par 
We will show that $V_i$, $V_8'$ 
are as described above. Denote by $C$ the image on $V_i$ 
(resp. $V_8'$) of the fibre of 
type $II$ (resp. $I_1$). Then one has 
$$-K_V \sim C \sim (9-i) \De, \,\, V = V_i, \, V_8'~ (i\neq 1)~ and~ -K_V\sim 4\De~ for~ i=1,$$
where $\De$ is the generator of the Weil divisor class group $Div(V)$.
The last assertion here comes from the observation that
the lattice on $X$ or $X'$ which is generated by
the section $E$, an elliptic fibre and the type $E_8$
sublattice in the type $II^*$ fibre, is unimodular
and hence equals $Div(X)$ or $Div(X')$.

By Kodaira's canonical bundle formula we have $K_X\sim -F$, 
where $F$ is any scheme-theoretic fiber of the elliptic fibration. 
This implies that $K_{V_i}\sim -C$ and hence $-K_{V_i}$ is ample. 
Similarly $-K_8'$ is ample. 

%
%
\par \vskip 1pc
When $i = 3$, the type $A_1, A_2$ singular points of $V_3$
lie on the smooth rational curve $\De$;
when $4 \le i \le 8$, each of $V_i$ and $V_8'$ 
has the unique singularity at
the cusp of the cuspidal rational curve $\De$.

We assert that $V_i \sms \De$ and $V_8' \sms \De$ are all
isomorphic to the affine plane ${\bf C}^2$.
Indeed, $V_i \setminus \De = V_8 \setminus \De$, 
and hence we have only to consider $V_8 \setminus \De$
($V_8' \setminus \De$ is similar).
Now $S_0 = 2(E + \sum_{i=1}^6 C_i) + C_3' + C_4'$
is the unique singular fibre of a ${\bf P}^1$-fibration 
$\varphi : X \rightarrow {\bf P}^1$ with $C_2'$ as a section. 
The assertion follows from the observation that 
$V_8 \setminus \De = X \setminus (S_0 + C_2')\cong \comx^2$.

\nind
{\bf Lemma 4.} (1) {\it $V_8$ and $V_8'$ are not isomorphic to each
other.}

\nind
(2) {\it Every rank $1$, Gorenstein log
del Pezzo surface $V$ satisfying $1\leq K_V^2\leq 6$ with $\pi_1(V^0) =
(1)$
is isomorphic to one of the $7$ surfaces $V_i~(3\leq i\leq 8),~ V_8'$.}

\nind
{\it Proof.} Some of the results and arguments 
in the proof below are well-known to the
experts, but we are giving them for the sake of completeness. 
At any rate, the assertion in part (2) above seems to be new.

Let $U \rightarrow V$ be a minimal resolution of singularities.

\nind
{\bf Claim 1.} $|-K_U|$ has a reduced irreducible member. 
Here we do not need the assumptions that $\rho(V) = 1$ and
$\pi_1(V^0) = (1)$.

\nind
Recall that by Kawamata-Viehweg vanishing theorem the group
$H^1(U,2K_U)=(0)$ [KMM, Theorem 1-2-3]. Hence by the Riemann-Roch
theorem,
one has dim $|-K_U| = K_U^2$. Suppose that a member $A$ of $|-K_U|$
contains an arithmetic genus $\ge 1$ irreducible component $A_0.$
The Riemann-Roch theorem implies that $|A_0 + K_U| \ne \emptyset$. From
this and $0 = A + K_U = (A_0 + K_U) + (A-A_0)$, we deduce that
$A_0 = A$ with $p_a(A_0) = 1$ and Claim 1 is true. 

\par
So we may assume that every member
of $|-K_U|$ is a union of smooth rational curves.
The Stein factorization and
the fact that $q(U) = 0$ imply that a general member of
$|-K_U|$ is of the form $M_1 + \cdots + M_k + F$,
where $F$ is the fixed part of the linear system, $M_i \cong {\bf P}^1$,
and $M_i \sim M_j$.

\par
Suppose that $K_U^2 = 1$. If $F = 0$, then $k = 1$
and $M_1^2 = 1$ and Claim 1 is true. 
Since $-K_U$ is nef and big, it is 1-connected by 
a result of C.P.Ramanujam.
Hence $1 = K_U^2 \ge (kM_1 + F) . kM_1 \ge  1 + k^2M_1^2$.
Thus $M_1^2 = 0, k = 1, M_1 . F = 1, K_U . F = 0$.
Now intersecting the relation $-K_U \sim M_1 + F$
with the smooth rational curve $M_1$ of self intersection 0,
one gets a contradiction. So Claim 1 is true when $K_U^2 = 1$.

\par
For $K_U^2 \ge 2$, let $U_1 \rightarrow U$ be
the blow-up of a point on $M_1 \setminus (M_2 + \cdots + M_k + F)$.
Then $-K_{U_1}$ is linearly equivalent to
the proper transform of $M_1 + \cdots + M_k + F$
and hence nef and big. If $U_1 \rightarrow V_1$ is the
contraction of all $(-2)$-curves then $V_1$ is a
Gorenstein log del Pezzo surface. For $U_1$,
we argue as in the case of $U$ and we can reduce
the proof of Claim 1 to the case $K_U^2 = 1$,
which has been dealt with in the previous paragraph.\\
This proves Claim 1.

\nind
We continue the proof of Lemma 4.
Suppose again first that $K_U^2 = 1$. Then dim $|-K_U|=1$ and hence $|-K_U|$
has a single base point, say $p.$ Let $A_1, A_2$ be two general members
of $|-K_U|$ meeting at $p$. Let $Y \rightarrow U$
be the blow-up of $p$ with $E$ the exceptional curve. Then $Y$
has a relatively minimal elliptic fibration $Y \rightarrow {\bf P}^1$
with the proper transforms of $A_1, A_2$ as fibres and $E$ as a section.

\nind
{\bf Claim 2.} (1) The singular fibre type of $Y \rightarrow {\bf P}^1$
is $II^* + II$ or $II^* + I_1 + I_1$. Hence $Y = X$ or
$Y = X'$ as described earlier in this section.

\nind
(2) The section $E$ is the only $(-1)$-curve on $Y$. 
All $(-2)$-curves are in the type $II^*$ fibre.
There are no other negative curves on $Y$.

\nind
By the assumption of Lemma 4, $Pic~ V$ is of rank 1. 
Since $V$ is simply connected, $Pic~V$ is also torsion free. 
Since $K_V^2 = 1$, one has $Pic~ V = {\zee} K_V$.
The assumption that $\pi_1(V^0) = (1)$ implies that the
Weil divisor class group $Div~(V)$ is torsion free and of rank 1
so that $Div~(V) = {\zee} C$ for some divisor $C$.
Write $C = a K_V$ with a rational number $a \leq 1$. Then
$a = C . K_V$ is an integer. So $a = 1$ and $Div~(V) = Pic~V$.
Note that $Div~(U)$ is the direct sum of the pull back of
$Div~(V)$ and the lattice generated by components of the 
exceptional divisor of the resolution $U \rightarrow V$. 
Now the unimodularity of $Div~(U)$ implies that $V$ has 
exactly one singularity and it is of type $E_8$.

\par
Clearly, the fibre on $Y$ containing the inverse of
the type $E_8$ divisor on $U$ (contractible to the singular
point on $V$) is of type $II^*$. There are no other 
reducible fibres by noting that $\rho(Y) = 10$ and that
the section $E$, a general fibre and the 8
components in the type $II^*$ fibre which is of Dynkin type $E_8$,
already give rise to 10 linearly independent classes of
$Div~(Y)$. The fact that the Euler
number of $Y$ is $12$ implies that the singular fibre type
of the elliptic fibration $Y \rightarrow {\bf P}^1$
is $II^* + II$ or $II^* + I_1 + I_1$.\\
This proves part (1) of Lemma 4 and part (1) of Claim 2.

Note that this also proves that $V_8,~V_8'$ are the only rank 1, 
Gorenstein log del Pezzo surfaces with $K_V^2=1.$

\par
Let $E_1$ be another $(-1)$-curve on $Y$. Then 
the observation that $-K_Y . E_1 = 1$ and the fact that 
$-K_Y$ is linearly equivalent to an elliptic fibre $F$
by Kodaira's canonical bundle formula imply that $E_1$
is a section of the elliptic fibration.
Hence we can write $E_1 = E + aF + D$ where
$a$ is rational and $D$ supported by the type $E_8$
divisor in the type $II^*$ fibre. Since $D^2 < 0$
and $D \cap E = \phi$ when $D \ne 0$, 
one has $D = 0$ by using $D$ to
intersect the expression of $E_1$. This leads to
$-1 = (E+aF)^2 = -1 + 2a$ and $a = 0$, a contradiction.

\nind
This proves (2) of Claim 2.

In view of what has been proved so far, 
we will assume that $2\leq d\leq 6.$ Denote $9-d$ by $c.$ From $d\leq 6$ we get $c\geq 3.$

\nind
{\bf Claim 3.} There is a composition of blow-ups $U_7 \rightarrow
U_6\rightarrow \cdots \rightarrow U_c = U$, so that if
$U_i \rightarrow W_i$ ($c\leq i \leq 7$) is the contraction 
of all $(-2)$-curves then $W_i$ is a rank 1, Gorenstein log del Pezzo
surface with $\pi_1(W_i^0) = (1)$ and $K_{W_i}^2 = 9 - i$.

\nind
Let $A$ be an irreducible member of $|-K_U|$ and $E_c$ a
$(-1)$-curve on $U$. Such a $(-1)$-curve exists because 
$K_V^2<7$ and hence $U$ is not a relatively minimal surface. 
Then $A$ meets $E_c$ at a point
$q$. Let $U_{c+1} \rightarrow U$ be the blow-up of
$q$ with $E_{c+1}$ the exceptional curve.
Then $-K_{U_{c+1}}$ is linearly equivalent to the
proper transform $A_{c+1}$ of $A$. Since $K_{c+1}^2>0,$ 
the divisor $-K_{c+1}$ is nef and big.
Moreover, the curves having 0 intersection with 
$A_{c+1}$ are precisely the inverse images of the $(-2)$-curves
on $U$ (contractible to singular points on $V$) 
and the proper transform $E_c'$ of $E_c$. It follows that 
the contraction of all the $(-2)$-curves on $U_{c+1}$ gives 
a surface $W_{c+1}$ as in the first part of Claim 3.

\nind
We will prove that $\pi_1(W_{c+1}^0) = (1)$.\\
Let $\Ga$ denote the union of the $(-2)$-curves in $U_c.$ 
Then $\Ga$ can be considered as a divisor on $U_{c+1}$ 
since $A$ is disjoint from $\Ga.$ By assumption, 
$\pi_1(U_c-\Ga)=\pi_1(U_{c+1}-\Ga)$ is trivial. 
Consider the natural map $\pi_1(U_{c+1}-\Ga-E_c')
\rar \pi_1(U_{c+1}-\Ga)=(1).$ By an application 
of Van Kampen theorem, we see that the kernel of this map 
is the normal subgroup of $\pi_1(U_{c+1}-\Ga-E_c')$ 
generated by a small loop around $E_c'.$ Since $E_{c+1}$ 
intersects $E_c'$ transversally once, this loop can be taken 
to be in $E_{c+1}.$ But $E_{c+1}-(\Ga \cup E_c')\cong {\comx}.$ 
Hence this loop is trivial in  $\pi_1(U_{c+1}-\Ga-E_c')$. 
This proves that $\pi_1(U_{c+1}-\Ga-E_c')$ is trivial.

\nind
We will now use the assertions of Claims 2 and 3  
to complete the proof of Lemma 4.

We will first prove that any rank 1, Gorenstein log del 
Pezzo surface $V$ with simply-connected smooth locus and 
$K_V^2=2$ is isomorphic to the surface $V_7$ constructed 
earlier in this section. We will show that $|-K_U|$ contains 
a cuspidal rational curve $C$ (of arithmetic genus $1$). 
Again let $A$ be an irreducible member of $|-K_U|$ and 
$L$ a $(-1)$-curve on $U.$ Let $U_8$ be the blow-up of $A\cap L.$ 
The contraction of all the $(-2)$-curves in $U_8$ produces 
either $V_8$ or $V_8',$ since these are the only surfaces 
with $d=1.$ If $U_8$ is the minimal resolution of $V_8$ 
then we already know that $|-K_{U_8}|$ contains a cuspidal curve. 
In fact, we know in this case by Claim 2 that there is a unique 
$(-1)$-curve in $U_8$ and the contraction of this curve is 
the minimal resolution of $V_7.$ In this case $V\cong V_7.$ 
So assume that $U_8$ is the minimal resolution of $V_8'.$ 
Recall that $X'$ is obtained by resolving the base locus of 
a pencil in $\proj^2$ generated by a nodal cubic $B'$ and $3$ times a 
line tangent at an inflexion point of the cubic.  

By the proof of Claim 2, the first $7$ blow-downs of $(-1)$-curves 
starting from $E$ are unique, viz. the curves $E,C_1,\ldots,C_6$ 
in this order. The contraction of $E$ produces $U_8.$ Hence the morphism 
$X'\rar \proj^2$ factors as $X' \rar U \rar\proj^2$. 
The existence of a cuspidal curve $C\in |-K_U|$ is 
equivalent to the existence of a cuspidal cubic $B$ in 
$\proj^2$ which has the same inflexion point and local intersection 
number $7$ at the inflection 
point with $B'.$ We can assume that the equation of $B'$ 
is $\{Y^2Z=X^3+XZ^2+\sqrt{-4/27}Z^3\},$ the point $[0,1,0]$ 
as the inflexion point, the tangent line being $\{Z=0\}.$ 
We can then take $B$ to be the cubic $\{Y^2Z=X^3\}.$ 
Now we see immediately that the blow-up of the point $B\cap L$ 
is the minimal resolution of $V_8.$ Then again $V$ is the surface $V_7.$ 

Now by Claim 2 we easily deduce that any rank 1, 
Gorenstein log del Pezzo surface with 
simply-connected smooth locus and $3\leq K^2\leq 7$ 
is one of the surfaces $V_i,~3\leq i\leq 7.$

\nind
Consider the degree $6$ hypersurface
$Z_a = \{W^2 + Z^3 + X^5Y + a X^4Z = 0\}$ in the weighted projective
space ${\proj}(1,1,2,3)$ with coordinates
$X, Y, Z, W$ of weights $1, 1, 2, 3$ respectively.

\nind
{\bf Lemma 5.} {\it When $a = 0$ (resp. $a \ne 0$),
$Z_a$ is isomorphic to $V_8$ (resp. $V_8'$).}

\nind
{\it Proof.} The affine open subset $\{X \ne 0\}$ of $Z_a$ is 
isomorphic to ${\comx}^2$. $Z_a$ has a unique singularity of 
type $E_8$ at $[0,1,0,0]$ (see [6]). It is easy to see that 
the boundary curve $\{X = 0 = W^2 + Z^3\}$ is isomorphic to 
a cuspidal cubic in ${\bf P}^2$; in particular it is irreducible 
and hence $Z_a$ has rank 1. If $a=0$ then the curve $C:=\{Y=0\}$ 
is a cuspidal rational curve and does not pass through the singular 
point of $Z_0.$ Futher, we have $K_{Z_0}\sim -C.$ 
Hence $Z_0$ is isomorphic to $V_8.$ If $a\neq 0$ 
then $C:=\{Y=0\}$ does not pass through the singular point 
of $Z_a.$ The curve $C=\{W^2+Z^3+aX^4Z=0\}$ in $\proj(1,2,3)$ 
is easily seen to be a smooth elliptic curve and $K_{Z_a}\sim -C.$ 
In [6] it is shown that the local rings of $Z_0$ and $Z_a$
at their singular points are not isomorphic. 
Hence $Z_a$ is isomorphic to $V_8'.$\\
This completes the proof of lemma 5.\\ 

We have also proved parts (i), (ii) of Theorem 2. The part (iii) is shown
in the construction of $V_i$, noting that $d = K_{V_i}^2 = 9-i$.

\nind
{\bf Proof of part (iv)}\\

Recall that $C$ is a cuspidal rational curve on $V_i$ not passing 
through the singular point of $V_i,$  where $3\leq i\leq 8$ 
and $-K_{V_i}\sim C.$ Then $C^2=9-i.$ By blowing up $V_i$ 
minimally at the singular point of $C$ we get a normal crossing divisor
 with smooth rational irreducible components $\cup_0^3B_i$ 
on a normal projective surface $V_i'',$ 
where $B_0^2=-1,B_1^2=-2,B_2^2=-3,B_3^2=3-i$, $B_0$ 
intersects $B_1,B_2,B_3$. The curves $B_1,B_2,B_3$ are mutually 
disjoint and $B_3$ is the proper transform 
of $C.$ Mumford's presentation for the fundamental group 
$G$ of the boundary of a nice tubular neighborhood of $\cup B_i$ 
is as follows (see [10]).\\
$$<e_2,e_3|~(e_2e_3)^2=e_2^3=e_3^{i-3}>$$

Now assume that $4\leq i\leq 8.$ Then this group is finite. 
This is the fundamental group at infinity of the affine surface $V_i-C.$ \\

\nind
{\bf Lemma 6.} {\it For $4\leq i\leq 8$ the surface $V_i-C$ 
is isomorphic to ${\comx}^2/G.$}\\

\nind
{\bf Proof.} Denote the surface $V_i-C$ by $S.$ We have proved 
above that $S$ has a unique singular point, say $p.$\\
First we treat the case $d=1.$ As seen above, in this case $S$ 
is isomorphic to the affine subset $\{W^2+Z^3+X^5=0\}$ of the 
projective surface $\{W^2+Z^3+X^5Y=0\}$ considered above given 
by $\{Y\neq 0\}.$ It is a classical fact that this $S$ is isomorphic
to the quotient of $\comx^2$ modulo the binary icosahedral 
group of order $120.$ Since $d=1,$ the fundamental group at 
infinity of $S$ has the presentation 
$<e_2,e_3|(e_2e_3)^2=e_2^2=e_3^5>.$ 
This group is the binary icosahedral group.\\

Now we assume that $4\leq i\leq 7.$\\

\nind
{\bf Claim.} For a small nice neighborhood $U$ of $p$ in $V_i,$ 
the natural homomorphism $\pi_1(U-\{p\})\rar \pi_1(S-\{p\})$ is 
an isomorphism.\\   
{\it Proof.} We will illustrate the proof by giving the argument 
for $d=2.$ Other cases are dealt in exactly the same way.\\
The surface $V_7$ is obtained from the elliptic surface $X$ 
by contracting the curves $E,C_1$ and all the other 
irreducible components of the type $II^*$ fiber $F_0$  
except for $C_2,$ giving rise to an $E_7$-singularity. 
Let $D$ denote the inverse image 
of $p$ in $X.$  Clearly, $S-\{p\}=X-(C\cup E\cup C_1\cup D).$
Let $N$ denote a suitable tubular neighborhood of $F_0$ in $X.$ 
It is easy to see that $N-(E\cup C_1\cup D)$ is a 
strong deformation retract of $X-(C\cup E\cup C_1\cup D).$ 
The neighborhood $N$
is a union of tubular neighborhoods $N_1,~N_2,~N_D$ of $ C_1,~C_2,~D$ 
respectively. Since $E\cap C_1$ is a single point 
$(N_2\cup N_D)-(C_1\cup D)$ is a strong deformation retract 
of $N-(E\cup C_1\cup D).$ Since $C_1\cap C_2$ is a single point $N_D-D$ is
a strong deformation retract of 
$(N_2\cup N_D)-(C_1\cup D).$ But $N_D-D$ is nothing but $U-\{p\}.$ \\
This proves the claim.\\

Let $W'$ denote the universal covering space of $S-\{p\}$ 
and let $W$ be the normalization of $S$ in the function 
field of $W'.$ By the claim just proved, $W$ contains a 
unique point, say $q,$ over $p$. This point is smooth by 
the claim just proved. Since 
$V_7$ has rank $1$ we see that $\chi(S-\{p\})=0$ where $\chi$ denotes the
topological Euler number. 
Therefore, $\chi(W')=0$ and hence $\chi(W)=1.$ Now $W$ 
is smooth, simply-connected and $b_2(W)=0,$ hence it is contractible.\\

Using $-K_{V_7}\sim C$ we see easily that the 
canonical bundle of $V_7''$ is linearly equivalent 
to a strictly negative linear combination of the curves 
$B_0,B_1,B_2,B_3.$ It follows that $\ka(S-\{p\})=-\infty.$ 
The map $W'\rar S-\{p\}$ being unramified
and proper we get $\ka(W')=-\infty.$ This implies that 
$\ka (W)=-\infty.$ Since $W$ is contractible, by 
a fundamental result of Miyanishi-Sugie-Fujita, $W$ is 
isomorphic to ${\comx}^2.$ This easily implies that 
$S\cong {\comx}^2/G,$ as desired.

\nind
This completes the proof of Theorem 2.

\section{Proof of Theorem 1}

Let $V$ be a Gorenstein normal projective surface such that there is a
surjective morphism $f: \proj^2 \rar V$.  We claim that $\pi_1(\vnot)$
is finite, where $\vnot = V \sms Sing~V$. Since $f^{-1}(Sing~V)$ has
codimension $2$ in $\proj^2,$ the complement $\proj^2 \sms
f^{-1}(Sing~V)$ is simply-connected.  If $W'$ is the universal cover
of $ \vnot$, then by the standard properties of topological coverings
the restriction of $f$ factors through $\proj^2 \sms f^{-1}(Sing~V)\rar W'$.
This easily proves the claim. Let $W \rar V$ be the
quasi-universal cover of $V$. Then $f$ factors as $\proj^2 \rar W \rar
V$. Now $\pi_1(W^\circ) = (1)$ and $W$ is Gorenstein of rank 1.  From
now onwards we will assume that $\vnot$ is simply-connected. 
 
\noin 
{\em First we deal with the case $d=1$}. \\

Assume that $V=V_8$. We will show that there is no non-constant morphism 
$\proj^2\rar V_8.$\\ 
Consider the surface $V_4.$ Let $p$ be the singular point and 
$C$ the cuspidal rational curve on $V_4^0$ such that 
$K_{V_4}\sim -C.$ As in Theorem 2, one has
$C \sim 5 \De$, where $\De$ is the image of the curve 
$C_5$ in $X$. By Theorem 2 we have a quotient map $\alpha$ of degree 
$5,~\comx^2 \rar V_4-C.$ Since $\alpha$ is proper, a neighborhood of
infinity of $\comx^2$ maps onto a neighborhood of infinity of $V_4 - C$.
Since $\comx^2$ is simply connected at infinity, we can use an argument
similar to the one at the beginning of this section to show that the order
of the fundamental group at infinity of $V_4 - C$ is at most 5.
There is a unique
point, say $\tilde p,$ 
in $\comx^2$ over $p.$ Let $Z$ be the normalization of 
$V_4$ in the function field of $\comx^2.$ We will show 
that $Z$ is isomorphic to $V_8.$ 

Since $V_4^0$ is simply-connected, the curve $C$ is ramified. 
Hence its inverse image in $Z$, say $D$, is irreducible 
and maps homeomorphically onto $C.$ If $q$ is the singular 
point of $C$ then the local analytic equation of $C$ at $q$ 
is $z_1^2+z_2^3=0.$ Hence the equation of $Z$ 
at its point $\tilde q$ over $
q$ is $z_1^2+z_2^3+z_3^5=0.$ This is the only singular 
point of $Z$ and it is clearly an $E_8$-singularity. 
As $Z$ contains $\comx^2,$ we see that $Z$ is a rank $1$, 
Gorenstein log del Pezzo surface such that $\pi_1(Z^0)=(1).$ 
To see that $Z$ is isomorphic
to $V_8$ we argue as follows.\\
There is a unique point in $Z$ over $C\cap \De.$ 
Since $\De-\{p\}-(C\cap\De)\cong \comx^*,$ 
we can prove easily that the inverse image of $\De$ 
in $Z,$ say $\tilde{\De},$ is irreducible and rational 
and it is smooth outside $\tilde p.$ Note that 
on $V_4$ we have $-K_{V_4} \sim 5\De\sim C.$ From the 
simple-connectedness of $Z^0$ we deduce that on $Z$ 
we have $D\sim \tilde{\De}$. This implies that on $Z$ 
there is a cuspidal rational curve $\tilde{\De}$ 
contained in the smooth locus of $Z$ such that $-K_Z \sim 
\tilde{\De}.$ By Theorem 2 we infer that $Z$ is 
isomorphic to $V_8.$ As remarked above we will 
show that there is no non-constant morphism $\proj^2\rar V_4.$ 
This will prove that $V_8$ is not dominated by $\proj^2.$\\ 

\nind
{\bf The cases $d=2,3,4,5.$}\\

First consider the cases $d = 3,4,5.$ By Lemma 2, 
$|-K_V|$ embeds $V$ in $\proj^d$ and the image is
projectively normal. Let $T \sst \comx^{d+1}$ be 
the cone over this embedding. The map 
$T \sms \{vertex\} \stackrel{\si}{\rar}V$ is a
locally trivial $\cstar$-bundle. It is the associated principal bundle
of the line bundle $\mp K_V$ over $V$. 
The sign does not play any significant role in 
our argument so we use the negative sign. This gives the exact sequence
$$ \pi_1(\cstar) \rar \pi_1(\tnot) \rar \pi_1(\vnot) \rar (1) $$
where $\tnot = T \sms Sing~T$. Since $\pi_1(\vnot) = 1$, we see that
$\pi_1(\cstar) \rar \pi_1(\tnot)$ is surjective.

\noin {\bf Lemma 7.} {\em The fundamental group of $\tnot$ is
isomorphic to $\zee/(d)$.} 

\noin {\bf Proof.} Let $C \sst \vnot$ be the cuspidal rational
curve with $C \sim -K_V$. 
Then $C^2 = d$. The inverse image $\si^{-1}(C)$ is the total
space of the principal bundle over $C$ of the line bundle $-K_V|_C$
which has degree $d$.  If $\olin{C} \stackrel{\beta}{\rar} C$ is the
normalization then $\beta$ is a homeomorphism. The pull-back of the
$\cstar$-bundle on $C$ to $\olin{C}$ is a $\cstar$-bundle of degree
$d$. Since $\olin{\comx} \cong \proj^1$, this pull-back has
fundamental group isomorphic to $\zee_d$.  It follows that
$\pi_1(T^\circ)$ is a homomorphic image of $\zee_d$.

We will now construct a cyclic $d$-fold \'{e}tale cover of
$T^\circ$. Since $\pi_1(\vnot) = (1)$ and rank of $V$ is one, the
class group of $V$ is cyclic. Let $\De$ be the generator 
of the Weil divisor class group as given in Theorem 2, so that 
$-K_V \sim d \De$ in $Div~(V)$, 
or ${\cal O}(-K_V) \cong {\cal O}(\De)^{\otimes d}$ restricted to
$V^\circ$.
We can find a suitable covering $\{U_i\}$ of $\vnot$ of Euclidean
balls and transition functions $f_{ij}$ for ${\cal O}(\De)|_{\vnot}$
such that $f_{ij}^d$ are the transition functions for
${\cal O}(-K_V)|_{\vnot}$. If ${\cal T} \rar \vnot$ is the total space of the
associated $\cstar$-bundle for ${\cal O}(\De)$, then the maps $U_i
\times \cstar \rar U_i \times \cstar$ defined by $(z,\la) \mapsto (z,
\la^d)$ patch to give an \'{e}tale cover of $T \sms Sing~T$ of degree
$d$. This proves the lemma.

The map $\proj^2 \rar V \sst \proj^d$ gives rise to a finite morphism
$\comx^3 \stackrel{\pi}{\rar} T$.  If $W \rar T$ is the
quasi-universal ${\zee}/(d)$-cover then $\pi$ factors as $\comx^3 \rar W
\stackrel{\tau}{\rar} T$. By lemma 3, $W$ admits an action of 
$\cstar.$ This is easily seen to be a
good $\cstar$-action such that $\comx^3 \rar W$ is a
$\cstar$-equivariant map. For a general fiber $F$ of $\si$ in
$T^\circ, \; \tau^{-1}(F)$ is smooth and irreducible because
$\pi_1(\cstar) \rar \pi_1(T^\circ)$ is surjective. Further, $(W \sms
\{vertex\}) // \cstar \cong V$.

We claim that the coordinate ring $\Ga(W)$ of $W$ is a UFD. 
Since the map $\comx^3\rar W$ is proper we see that $Div~(W)$ 
is finite. Any non-trivial torsion line bundle on $W^0$ 
gives rise to a non-trivial topological covering of $W^0.$ 
But $\pi_1(W \sms Sing~W) = (1).$ Hence $\Ga(W)$ is a UFD. Recall that
$d = 3,4$ or $5.$ Then $G/\Ga \cong \zee/(d)$, where $\Ga = [G,G]$
and $G$ is defined in Lemma 6.

Let $Z$ be the inverse image of $V \sms C$ in $W \sms
\{vertex\}$. Then $Z \sms l$ is the total space of the principal
$\cstar$-bundle associated to the line bundle 
${\cal O}(\De)|_{V \setminus C}$, 
where $l$ is the inverse image in $Z$ of the singular point of $V$.
Consider the following action of $G$ on $\comx^2 \times \cstar$. 
We consider the surjection $\mu : G \rar G/\Ga (\cong \zee/(d))$. 
Let $\olin{g}$ be a generator of $G/\Ga$ and $g$ a 
lift of $\olin{g}$ in $G$. Any element of $G$ has the unique
expression $\ga.g^b, \; \ga \in \Ga$ and $0 \leq b \leq d-1$. 
Let $h = \ga g^b$ act on $(z, \la)$ by 
$h(z, \la) = (hz, \omega^b \la)$, 
where $\omega=exp(2\pi i/d).$ 

\noin {\bf Lemma 8.} {\em With the above action we have 
$(\comx^2 \times \cstar)/G \cong Z$.}  

\noin {\bf Proof.} Recall that the inverse image of $V \sms C$ in $T$
is isomorphic to $(V \sms C) \times \cstar \approx (\comx^2/G) \times
\cstar$. This is because $K_V|_{V \sms C}$ is a trivial line
bundle. Define the map $\alpha : \comx^2 \times \cstar \rar
(\comx^2/G) \times \cstar$ given by $(z, \la) \mapsto (\olin{z},
\la^d)$. If $\ga g^b$ is an arbitrary element of $G$ as above, then
$\alpha ((\ga g^b)(z, \la)) =
(\olin{z}, \la^d)$. Hence the map $\alpha$ factors through 
$\olin{\alpha} : (\comx^2 \times \cstar)/G \rar (\comx^2/G) \times \cstar$. 
If $q \in \comx^2/G$ is a
smooth point of $\comx^2/G$ then $\alpha^{-1}(q,\la)$ has $d |G|$
points. Hence $\olin{\alpha}^{-1}(q,\la)$ has $d$ 
points and $\olin{\alpha}$ is \'{e}tale outside 
$\{\olin{0}\} \times \cstar$. We
see easily that the inverse image of an orbit $q \times \cstar$ in
$(\comx^2 \times \cstar)/G$ is connected. On the other hand, for the
map $W \sms \{vertex\} \rar T \sms \{vertex\}$, we have proved that
any good orbit $q \times \cstar \sst T \sms \{vertex\}$ lifts to a
single orbit in $W \sms \{vertex\}$.  From these two observations, we
deduce that $\comx^2 \times \cstar/G$ is naturally 
isomorphic to $Z$, proving the lemma.

We now come to the {\em punch line}
(deducing a contradiction to the fact that
the affine open subset $Z \subseteq W$
also has a UFD as the coordinate ring).

\noin {\bf Lemma 9.} {\em The affine 3-fold 
$(\comx^2 \times \cstar)/G$ is not a UFD.}

\noin {\bf Proof:} The map 
$\comx^2 \times \cstar \rar (\comx^2 \times
\cstar)/G$ is unramified outside $\{\olin{0}\} \times \cstar$, 
which has codimension $2$.  Let
$U$ be the group of units in the coordinate ring $R$ of $\comx^2
\times \cstar$. Then $U \cong \cstar \times \zee t$. Here $\cstar$ is
the multiplicative unit group of the underlying base field 
and $\zee t$ generates the group of units of $\cstar$ (modulo the
non-zero constants). By Samuel's
descent
theory 
(see [12], Chapter III), the divisor
class group of $(\comx^2 \times \cstar)/G$ is $H^1(G,U)$. Consider the
short exact sequence of $G$-modules
$$(1)\rar {\comx}^*\rar U\rar {\zee} t\rar (1)$$
Here, we consider $U/\cstar\cong {\zee}t$ as a $G$-module. The long exact 
cohomology sequence corresponding to this looks like
$$(1)\rar H^0(G,{\comx}^*)\rar H^0(G,U)\rar H^0(G,{\zee} t)\rar
H^1(G,{\comx}^*)\rar H^1(G,U)\rar\cdots$$
Now $H^0(G,{\zee}t)$ is the cyclic subgroup of ${\zee}t$ of index $d$ 
invariant under $G$ and $H^0(G,U)$ is a direct sum 
$H^0(G,\cstar)\oplus H^0(G,{\zee}t)$. Hence $H^1(G, \cstar)$ is a 
subgroup of $H^1(G,U)$. On the other hand, since $G$ acts trivially 
on the field of constants $\comx$, $H^1(G,\cstar) =
Hom(G,\cstar) \cong Hom(G/\Ga, \cstar) \cong \zee/(d)$. But $d
\neq 1$ by assumption. This proves the lemma.

So we have proved that $V$ can not be an image of $\proj^2$,
for those $V = V_i$ with $d = K_{V_i}^2 = 9 - i = 1, 3, 4, 5$.

Next we consider the cases $d = 2, 6, 8$. 

\noin {\em The case $d = 2$}. \\
In this case $|-K_V|$ does not give an embedding of $V$ but $|-2K_V|$
gives a projectively normal embedding of $V$ in $\proj^6$. Let $T$
denote the affine cone over this embedding. Then we work with this
cone $T$. An easy modification of the argument for $3 \leq d \leq 5$
shows that $V_i$ ($i = 9 - d = 7$) is not an image of $\proj^2$.

\noin {\em The case $d = 6$}. \\
In this case $V \sms C$ is not isomorphic to $\comx^2/G$. In fact $V
\sms C$ has singularity type $A_1+A_2$. Consider the action of
$\zee_6$ on $\proj^2$ given by $\si[x_0, x_1, x_2] = [x_0, \omega x_1,
-x_2]$ where $\omega = exp(2 \pi i/3)$. Let $x,y$ be
suitable affine coordinates on $\proj^2 \sms \{x_0 = 0\}$. Then
$\si(x,y) = (-x,\omega y)$ and $\si^2(x,y) = (x, \omega^2 y)$ which
implies that $\si^2$ is a pseudo-reflection. Similarly $\si^3$ is a
pseudo-reflection. Hence $\comx^2/\langle \si \rangle$ 
is smooth and hence isomorphic to $\comx^2$.

Consider $\proj^2 \sms \{x_1 = 0 \}$. Now $\si(x,y) = (\omega^2 x,
-\omega^2 y)$ and $\si^3(x,y) = (x,-y)$ and hence $\si^3$ is a
pseudo-reflection. The ring of invariants of $\si^3$ is $\comx[x,y^2]$. 
The action of $\olin{\si}$ on $\comx[x,y^2]$ is $(x,y^2) \mapsto
(\omega^2x, \omega y^2)$. This gives an $A_2$-singularity.

Finally consider $\proj^2 \sms \{x_2 = 0\}$. In this case $\si(x,y)
= (-x, -\omega y)$ which implies that $\si^2(x,y) =
(x,\omega^2y)$. The invariants are $\comx[x,y^3]$ and hence
$\olin{\si}(x,y^3) = (-x, -y^3)$. This gives an $A_1$-singularity on
the quotient.

By Theorem 2, $\proj^2/\langle \si \rangle = V_3$
where $d = K_{V_3}^2 = 9 - 3 = 6$. In other words,
$V_3$ ($d = 6$) is the quotient of $\proj^2$ by an action of $\zee/(6)$.

\noin {\em The case $d = 8$.} \\
In this case $V$ is the quadric cone $Q$ in $\proj^3$. Let $\zee/(2)$ act
on
$\proj^2$ by $\si[x_0, x_1, x_2] = [-x_0, -x_1, x_3]$. Then $V \cong
\proj^2 /\zee/(2)$.\\
In conclusion, we have so far proved that if $V$ is a rank 1, 
Gorenstein log del Pezzo surface such that $\pi_1(V^0)$ is 
trivial and there is a non-constant morphism $\proj^2\rar V$ 
then $V$ is isomorphic to either $\proj^2,~Q$, a surface of singularity
type $A_1+A_2$, or $V_8'$. The values of $d$ in these cases are $9,8,6,1$ 
respectively.\\
We will next prove that if $V$ is a rank 1, Gorenstein log del 
Pezzo surface dominated by $\proj^2$ such that $\pi_1(V^0)$ 
is non-trivial then its quasi-universal cover is either 
$\proj^2$ or $Q.$\\
Let $W$ be the quasi-universal cover of $V$ and $g:W\rar V$ 
the covering map. Then $W$ is a rank 1 Gorenstein, log del Pezzo surface
dominated by $\proj^2$ such that $\pi_1(W^\circ) = (1)$. \\
Assume first that $K_W^2=1.$ Since $K_W\sim g^*K_V$, we get 
$K_W^2=1=~deg~g\cdot K_V^2$. This means that $g$ is an isomorphism. 
Hence $V_8,V_8'$ cannot occur as the quasi-universal cover of $V$.\\
Suppose that $W$ is of singularity 
type $A_1 + A_2.$ Let the singular points of $W$ be $p,q.$ 
As $g$ is unramified over $V^\circ$, the images 
$p':=g(p),~q':=g(q)$ are singular points
of $V.$ Also, $6 = K_W^2=(deg~g) K_V^2.$ We analyse the possible cases.\\
{\it Case 1.} Suppose that $deg~g=2.$ 
Then any singular point of $V$ other than $p',q'$ is an $A_1$-singular
point. 
Let $\vtld\rar V$ be a minimal resolution of 
singularities. Since $K_V^2=3$ and $\rho(V) = 1$, the number of 
irreducible exceptional curves for the map $\vtld\rar V$ 
is $6.$ It is easy to see that $p',q'$ are of type $A_3,A_5$ 
respectively. This is a contradiction.\\
{\it Case 2.} Suppose that $deg~g=3.$ Then $K_V^2=2$ 
and the number of exceptional curves for $\vtld\rar V$ is $7.$ 
The type of $p'$ is $A_5$ and $q'$ will give rise to 
more than $2$ exceptional curves.  This is a contradiction.\\
{\it Case 3.} Suppose that $deg~g=6.$ Now $K_V^2=1$ 
and the number of exceptional curves is $8.$ 
The local fundamental groups at $p',q'$ 
have orders $12,~18$ respectively. Looking at the 
possible Dynkin types of these singularities we arrive at a contradiction.\\
This proves part (2) of Theorem 1.\\

\nind
{\bf Proof of Part (3).}\\

As above, let $f:\proj^2\rar V$ be a non-constant morphism.  We will 
assume that $V$ is not $V_8'$. Assume that the quasi-universal cover of $V$ is $Q.$ We will 
prove that $V$ is isomorphic to $\proj^2/H$ for a suitable 
finite group of automorphisms $H$ of $\proj^2.$ Of course this $\proj^2$
may not be the same as the projective plane dominating $V$. \\

Let $q$ be the singular point of $Q.$ It is easy to 
see that $Q$ contains a smooth rational curve $D$ with $D^2=2$ 
and not passing through $q.$ Further, $\pi_1(Q-D-\{q\})={\zee}/(2).$ 
If $Y$ is the universal cover of $Q-D-\{q\}$ then the normalization of 
$Q$ in the function field of $Y$ is isomorphic to 
$\proj^2$ such that the inverse image of $D$ in $\proj^2$ is a line.  
We will use this observation below.\\

By assumption, $V$ is isomorphic to $Q/G$ with $G = \pi_1(V^0)$,
and the map $Q\rar V$ is unramified over $V^0.$ Then $K_Q^2=8=|G|K_V^2.$ 
Therefore $|G|$ is of order $2,4$ or $8.$ The action of $G$ 
extends uniquely to the minimal resolution of singularities 
of $Q,$ viz. to Hirzebruch surface 
$\Sigma_2.$ Let $M$ be the unique $(-2)$-curve 
on $\Sigma_2$ and let $L$ be a fiber of the $\proj^1$-fibration 
on $\Sigma_2.$ Then $G$ acts naturally on $|D|$ and $D\sim M+2L.$ 
The linear system $|D|$ is parametrized by $\proj^3.$ The 
subspace of members of this of the form $M+2L',$ where $L'$ 
is a fiber of the $\proj^1$-fibration is parametrized by 
$\proj^2.$ This $2$-dimensional subspace is clearly stable 
under the action of $G.$ The complement of this $2$-dimensional 
subspace in $|D|$ is parametrized by $\comx^3$. 
As $G$ is a finite $2$-group, 
by a standard result in Smith Theory the action of $G$ on 
$\comx^3$ has a fixed point (see, [2]). This means that 
there is an irreducible smooth rational curve $D_0 \in |D|$ which 
is stable under $G.$ Hence the set 
$Q-D_0-\{q\}$ is also $G$-stable. This implies 
that the map $Q\rar V$ is unramified over $Q-D_0-\{q\} - A_2$
with a finite set $A_2$ and 
$Q-D_0-\{q\} - A_2$ is the full inverse image of its image in $V.$ 
>From the observation made in the beginning, we conclude that 
there is a line $C$ in $\proj^2$ such that the composite map 
$\proj^2-C-\{q_1\} - A_1 \rar V-E-\{q_3\} - A_3$ is unramified, 
where $E, q_3, A_3$ (resp. $C, q_1, A_1$) 
are images (resp. inverse images) of $D, q, A_2$ in $V$ (resp. $\proj^2$). 
Since $\proj^2-C-\{q_1\} - A_1$ is simply-connected, 
the morphism $\proj^2\rar V$ is a Galois map.\\

To complete the proof of part (3), we will now show that $V_8'$ is not 
isomorphic to a quotient $\proj^2/G$.\\
Suppose that $V=V_8'\cong \proj^2/G$ and let $f:\proj^2\rar V$ be the 
quotient map and $deg~f=n$. Let $\Gamma_1,\ldots,\Gamma_m$ be the 
irreducible components of the branch locus in $V.$ Denote by 
$\Gamma_{ij}$ the irreducible components of $f^{-1}\Gamma_
i$ with ramification index $e_i$. 
Then $f^*\Gamma_i=\Sigma_j e_{i}\Gamma_{ij}$. For the canonical bundle 
we have $K_{\proj^2}\sim f^*K_V+\Sigma_{i,j} (e_i-1)\Gamma_{ij}$. 
Write $\Gamma_i\sim \delta_i C$, where $C$ is an irreducible curve on $V$
such that $K_V\sim -C$. This gives $K_{\proj^2}\sim \Sigma_i
f^*({\frac{e_i-1}{e_i}}\delta_i-1)C$. Since $K_{\proj^2}$ is negative, we
infer easily that the branch curve in $V$ is irreducible and $\delta_i=1$.
This means that the branch curve is a member of $|-K_V|$. But it can be easily
seen from the arguments used earlier that for any member $D$ of $|-K_V|$ the
complement $V\sms D$ is simply-connected. This shows that $V_8'$ is not a
quotient of $\proj^2$.\\   
This also completes the proof of Theorem 1.\\

\nind
\section{Classification of Gorenstein Quotients of $\proj^2$}

Let $V=\proj^2/G$ be a normal Gorenstein quotient of $\proj^2.$\\

\nind
{\it Case 1.} The quotient map $f:\proj^2\rar V$ is a 
quasi-universal covering.\\
In this case $9 = (deg~f) K_V^2.$ Hence $|G|=3$ or $9.$\\
{\it Case 1.1.} Suppose that $|G|=3.$\\
In this case every singular point of $V$ is an $A_2$-singularity. 
Let $\vtld\rar V$ be the minimal resolution of singularities. 
Since $K_V^2=3,$ the number of irreducible components of the 
exceptional divisor for the map $\vtld\rar V$ is $6.$ This means
that $V$ has exactly three singular points.
It is easy to see that in a suitable coordinate system,
the action of a generator of $G$ is 
given by $[X,Y,Z]\rar [X,\omega Y,\omega^2 Z],$ 
where $\om$ is a primitive cubic root of unity.\\
{\it Case 1.2.} Suppose that $|G|=9.$\\
Now $K_V^2=1.$ The number of irreducible components of the 
exceptional divisor for the map $\vtld\rar V$ is $8.$ 
The map $f$ is unramified outside finitely many points. Hence the order of
the local fundamental group at any singular point of $V$ divides 9. First,
we claim that $A_8$ cannot occur as a singularity of $V$. For otherwise
$V$ has no other singularity and the map $f$ is unramified outside the
singular point and there is a unique singular point in $\proj^2$ over this
point. But $\chi (\proj^2 \sms \{one~point\})= 2 = \chi (V \sms \{
one~point\})$. This contradicts the multiplicativity of $\chi$ for
topological coverings. We deduce that $V$ has exactly four $A_2$-type 
singular 
points. There are exactly
three distinct points in $\proj^2$ over each of these.\\
We claim that $G$ is isomorpic to a direct sum $\zee/(3)\oplus \zee/(3)$
so that $\pi_1(V^0) \cong \zee/(3) \oplus \zee/(3)$.\\
Assume that $G$ is cyclic, say $G = \langle g \rangle.$ We can assume 
that the action of $g$ sends $[X,Y,Z] \rar [\omega X,\omega^q Y,Z]$ 
for some integer $q,$ where $\omega = exp(2 \pi \imath /9)$.
The point $[0,0,1]$ is fixed under 
the group and the action of $G$ near this point cannot 
have any non-trivial pseudo-reflection 
as the map $f$ is divisorially unramified. 
Since the singularities in $V$ are Gorenstein we deduce that 
$q=8.$ But this produces an $A_8$-type singularity. 
This proves the claim. 

Now we will construct an explicit example giving such a surface.
Let $G = \zee/(3) \oplus \zee/(3) = (g_1) \oplus (g_2)$, where $g_1$ sends
$[X, Y, Z] \mapsto [X, \omega Y, \omega^2 Z]$ where $\omega$ is a 
primitive cube root of unity and $g_2$ sends $[X, Y, Z] \mapsto [Z, X,
Y]$. Let $W = \proj^2/(g_1)$. The points [1, 0, 0], [0, 1, 0] and [0, 0,
1] are fixed by $g_1$ and form a single $g_2$-orbit. Their images in $W$
are $A_2$ type singular points. The points $[1,1,1]$, [1, $\omega$,
$\omega^2$] and [1, $\omega^2$, $\omega$] are fixed by $g_2$ and form a
single $g_1$-orbit. Their image in $W$ is another $A_2$ type singular
point. The points [1, 1, $\omega$], [1, $\omega$, 1] and [1, $\omega^2$,
$\omega^2$] form both a single $g_1$-orbit and a single $g_2$-orbit. Their
image in $V$ is a singular point of type $A_2$. Finally [1, $\om^2$, 1],
[1, 1, $\om^2$], [1, $\om$, $\om$] form both a single $g_1$-orbit and a
single $g_2$-orbit. Their image in $V$ is the fourth singular point of
type $A_2$ in $V$.

\nind
{\it Case 2.}  $h:Q\rar V$ is the quasi-universal covering. Let $H$ be the
Galois group. \\
Let $q$ be the singular point of $Q.$ 
In this case, $8 = (deg~g) K_V^2.$\\
{\it Case 2.1.} Suppose that $|H|=2.$\\
Since $K_V^2=4$ the number of exceptional irreducible components 
for $\vtld \rar V$ is $5.$ The image of $q$ is of $A_3$-type. 
All other singular points are of $A_1$-type. Thus $V$ 
has singularity type $A_3, A_1, A_1$.
By Theorem 1, $V \cong \proj^2/G$ with
$|G| = 4$. The existence of type $A_3$ singularity
on $V$ shows that $G \cong \zee/(4)$. \\
An explicit example in this case is given as follows.
Let $g$ be the automorphism of $\proj^2$ of order 4 sending $[X, Y, Z]
\mapsto [X, \imath Y, -\imath Z]$ where $\imath$ is a square root of $-1$.
The point $[1,0,0]$ is fixed by $g$ and every point on $X = 0$ is fixed
by $g^2$. No
other point of $\proj^2$ has a non-trivial isotropy group. The image of
$[1,0,0]$ on $V$ is an $A_3$ type singularity and the images of $[0,1,0]$
and
$[0,0,1]$ on $V$ are $A_1$ type singular points. The quotient
$\proj^2/(g^2)$ is the quadric $Q$. \\
{\it Case 2.2.} Suppose that $|H|=4.$\\
Now $K_V^2=2.$ The number of irreducible exceptional 
curves for $\vtld\rar V$ is $7.$ The image of $q$ in $V,$ 
say $q',$ has local fundamental group of order $8.$ 
If $q'$ is of type $A_7$ then $V$ has no other singular points. 
But this contradicts the multiplicativity of $\chi.$\\
Hence $q'$ is of $D_4$-type. Then the other singular points 
of $V$ are of type $A_1,A_2$ or $A_1,A_1,A_1$. But the order 
of the local fundamental group of any other singular point is 
a divisor of $4$. Hence we conclude that $V$ has singularity  type 
$D_4+~3A_1.$ By Theorem 1, $V \cong \proj^2/G$ with
$|G| = 8$. The existence of the type $D_4$ singularity
on $V$ shows that $G$ is the binary dihedral group of order 8. An explict
action of $G$ is as follows. Let $g_1$ map $[X, Y, Z] \mapsto [X, \imath 
Y, -\imath Z]$ and $g_2$ map $[X, Y, Z] \mapsto [X, \imath Z, \imath Y]$.
Then $g_1,~g_2$
generate a group of order 8 such that the subgroup $(g_1)$ has index 2 and
hence normal in $G$. Arguing as in case 1.2 by considering the
intermediate quotient $W = \proj^2/(g_1)$ we see that $V$ has singularity
type $D_4+~3A_1$. \\
{\it Case 2.3.} Suppose that $|G|=8.$\\
Now $K_V^2=1$ and the number of irreducible exceptional curves 
for $\vtld\rar V$ is $8.$ The image $q'$ of $q$ has local 
fundamental group of order $16.$ This cannot be of $A_{15}$-type. 
Hence $q'$ is of type $D_6.$ The other singular points are of type $
A_2$ or $A_1,A_1.$ Again $A_2$ cannot occur because the order of the 
local fundamental group is not a divisor of $8.$ Hence $V$ is of 
singularity type $D_6+~2A_1.$ We claim that this case cannot occur. 
Over a singular point of type $A_1$ there are four 
points in $Q.$ This easily contradicts the multiplicativity of $\chi.$\\

Finally, we consider the surface $V_3$ of singularity type $A_1~+~A_2$. An
explicit action of $\zee/(6)$ on $\proj^2$ which produces this quotient is
as follows. Let $g$ be an automorphism of $\proj^2$ of order 6 sending
$[X, Y, Z] \mapsto [X, -Y, \om Z]$, where $\om$ is a primitive cube root
of unity. The image of [0, 1, 0] in $V$ is an
$A_1$ type singularity. The image of [0,0,1] in $V$ is an $A_2$ type
singularity.

\nind
We have thus proved the following result.\\

\nind
{\bf Lemma 10.} {\it If ${\proj}^2$ is the quasi-universal cover 
of a normal, Gorenstein, projective surface $V$ (not isomorphic to
$\proj^2$) 
then either $V^0$ has fundamental group $\zee/(3)$ and singularity 
type $3A_2$, or the fundamental group of $V^0$ is 
$\zee/(3)\oplus\zee/(3)$ and $V$ is of singularity type $4A_2.$\\
If $Q$ is the quasi-universal cover of a normal 
projective surface $V$ (not isomorphic to $Q$) then either 
the order of the fundamental group of $V^0$ is $2$,
$V$ is of singularity type $A_3+2 A_1$ and 
$V \cong \proj^2/\zee/(4)$,
or the order of the fundamental group of $V^0$ is $4$,
$V$ is of singularity type $D_4+3A_1$
and $V \cong \proj^2/H$ with $H$ the quaternion
group of order $8$.\\
The only other Gorenstein surface not covered by the 
above cases which is isomorphic to a quotient of 
$\proj^2$ is the surface $V_3$ of singularity type $A_1+A_2$. 
The fundamental group of $V_3^0$ is trivial
and $V_3 \cong \proj^2/\zee/(6)$.}

\section{Log del Pezzo non-Gorenstein case}
Assume that $f:\proj^2 \rar V$ is surjective and $V$ is a log del Pezzo
surface such that $\pi_1(\vnot) = (1)$.  Let $V \sst \proj^N$ be a
suitable projectively normal embedding and $T$ the cone over $V$ in
$\comx^{N+1}$. Then we get a finite map $\comx^3 \rar T$. Denote by
$W$ the quasi-universal cover of $W$. As before $T$ admits 
a good ${\bf C}^*$-action such that the map $\comx^3 \rar T$ 
is $\cstar$-equivariant.

\noin {\em Conjecture 1:} Let $Y$ be a normal affine variety with a
good $\cstar$-action. Suppose $\pi: \comx^n \rar Y$ is a proper
surjective morphism which is $\cstar$-equivariant. Assume that
$\pi_1(Y \sms Sing~Y) = (1)$. Then $Y$ is isomorphic to $\comx^n$ with
a suitable good $\cstar$-action.

If this conjecture has an affirmative answer then $Y//\cstar$ is a
weighted projective space $\proj(a,b,c,\ldots)$. In the case of log
del Pezzo surface under consideration, we know that $W//\cstar \cong
V$. Hence $V$ is isomorphic to $\proj(a,b,c)$. Therefore an
affirmative answer to the above conjecture gives an affirmative answer
to

\noin {\em Conjecture 2:} Let $V$ be a log del Pezzo surface with a
surjective morphism $f : \proj^2 \rar V$. Then $V$ is isomorphic to a
quotient $\proj(a,b,c)/G$, with $G$ isomorphic to the fundamental
group $\pi_1(\vnot)$.\\
In particular, if $\vnot$ is simply-connected then $V$ is isomorphic to
$\proj(a,b,c)$.

\noin {\bf Remark.} In the Gorenstein case, $d = 8$ corresponds to
$\proj(1,1,2)$ and $d = 6$ to $\proj(1,2,3)$.

\noindent
R.V. Gurjar, School of Mathematics, 
Tata Institute of Fundamental research, 
Homi-Bhabha road, Mumbai 400005, India. (e-mail: 
gurjar@math.tifr.res.in)\\
C.R. Pradeep, School of Mathematics, 
Tata Institute of Fundamental research, 
Homi-Bhabha road, Mumbai 400005, India. (e-mail:
pradeep@math.tifr.res.in)\\
D.-Q. Zhang, Dept. of Mathematics, 
National University of Singapore, 2 Science Drive 2, Singapore 117543. 
(e-mail: matzdq@math.nus.edu.sg)\\
\end{small}

\end{document}